\documentclass{autart}    

\usepackage[dvips]{epsfig}     
\usepackage{amsmath}
\usepackage{amssymb}

\usepackage{amsfonts,bbold,graphicx,picins}
\usepackage{color,natbib,url}

\newtheorem{definition}{Definition}
\newtheorem{theorem}{Theorem}
\newcommand{\wno}{\mbox{wno}}

\begin{document}

\begin{frontmatter}
\title{A Data-driven Approach to Robust Control of Multivariable Systems by Convex Optimization}

\author{Alireza Karimi\thanksref{Corauthor}}
\author{and Christoph Kammer}

\thanks[Corauthor]{Corresponding author: alireza.karimi@epfl.ch}

\address{Laboratoire d'Automatique, \'Ecole Polytechnique F\'ed\'erale de Lausanne (EPFL), CH-1015 Lausanne (Switzerland)}


\begin{abstract}
The frequency-domain data of a multivariable system in different operating points is used to design a robust controller with respect to the measurement noise and multimodel uncertainty. The controller is fully parametrized in terms of matrix polynomial functions and can be formulated as a centralized, decentralized or distributed controller. All standard performance specifications like $H_2$, $H_\infty$ and loop shaping are considered in a unified framework for continuous- and discrete-time systems. The control problem is formulated as a convex-concave optimization problem and then convexified by linearization of the concave part around an initial controller. The performance criterion converges monotonically to a local optimal solution in an iterative algorithm.  The effectiveness of the method is compared with fixed-structure controller design methods based on non-smooth optimization via multiple simulation examples. 
\end{abstract}

\end{frontmatter}
%
\section{Introduction}\label{sec1}
Recent developments in the fields of numerical optimization, computer and sensor technology have led to a significant reduction of the computational time of optimization algorithms and have increased the availability of large amounts of measured data during a system's operation.  These progresses make computationally demanding data-driven control design approaches an interesting alternative to the classical model-based control problems. In these approaches, the controller parameters are directly computed by minimizing a control criterion which is a function of measured data. Therefore, a parametric model of the plant is not required and there are no unmodeled dynamics. The only source of uncertainty is the measurement noise, whose influence can be reduced significantly if the amount of measurement data is large.


Frequency-domain data is used in the classical loop-shaping methods for computing simple lead-lag or PID controllers for SISO stable plants. The Quantitative Feedback Theory (QFT) uses also the frequency response of the plant model to compute robust controllers \citep{Hor93}. In these approaches the controller parameters are tuned manually using graphical methods. New optimization-based algorithms have also been proposed recently \citep{MABH16}. The set of all stabilizing PID controllers with $H_\infty$ performance is obtained using only the frequency-domain data in \cite{KB08}.  This method is extended to design of fixed-order linearly parameterized controllers in \cite{PK15,PK16}. The frequency response data are used in \cite{HDAV10} to compute the frequency response of a controller that achieves a desired closed-loop pole location. A data-driven synthesis methodology for fixed structure  controller design problems with  $H_\infty$ performance is presented in  \cite{DWS09}. This method uses the $Q$ parameterization in the frequency domain and solves a non-convex optimization problem to find a local optimum. Another frequency-domain approach is presented in \cite{KNND13a} to design reduced order controllers with guaranteed bounded error on the difference between the desired and achieved magnitude of sensitivity functions. This approach also uses a non-convex optimization method. 

Another direction for robust controller design based on frequency-domain data is the use of convex optimization methods. A linear programming approach is used to compute linearly parametrized (LP) controllers for SISO systems with specifications in gain and phase margin as well as the desired closed-loop bandwidth in \cite{KKL07,Sae14}. A convex optimization approach is used to design LP controllers with loop shaping and $H_\infty$ performance in \cite{KG10}. This method is extended to MIMO systems for computing decoupling LP-MIMO controllers in \cite{GKL10b}. 
Recently, the necessary and sufficient conditions for the existence of data-driven $H_\infty$ controllers for SISO systems has been proposed in \cite{KNZ16}. 

The use of the frequency response for computing SISO-PID controllers by convex optimization is proposed in \cite{HAB13}. This method uses the same type of linearization of the constraints as in \cite{KG10} but interprets it as a convex-concave approximation technique. An extension of \cite{HAB13} for the design of  MIMO-PID controllers by linearization of quadratic matrix inequalities is proposed in \cite{BHA16} for stable plants. A similar approach, with the same type of linearization, is used in \cite{SOW10} for designing LP-MIMO controllers (which includes PID controllers as a special case). This approach is not limited to stable plants and includes the conditions for the stability of the closed-loop system.

In this paper, a new data-driven controller design approach is proposed  based on the frequency response of multivariable systems and convex optimization. Contrarily to the existing results in \cite{GKL10b,BHA16,SOW10}, the controller is fully parameterized and the design is not restricted to LP or PID controllers. The other contribution is that the control specification is not limited to $H_\infty$ performance. The $H_2$, $H_\infty$ and mixed $H_2 / H_\infty$ control problem as well as loop shaping in two- and infinity-norm are presented in a unified framework for systems with multimodel uncertainty. A new closed-loop stability proof based on the Nyquist stability criterion is also given. 

It should be mentioned that the problem is convexified using the same type of approximation as the one used in \cite{BHA16,SOW10}. Therefore, like other fixed-structure controller design methods (model-based or data-driven), the results are local and depend on the initialization of the algorithm.


\section{Preliminaries}\label{sec2}
The system to be controlled is a Linear Time-Invariant Multi-Input Multi-Output (LTI-MIMO) system represented by a multivariable frequency response model $G(e^{j\omega}) \in \mathbb{C}^{n \times m}$, where $n$ is the number of outputs and $m$ the number of inputs. The frequency response model can be identified using the Fourier analysis method from $m$ sets of input/output sampled data as \citep{PS01}:
\begin{equation}
G(e^{j\omega})=\left[ \sum_{k=0}^{N-1} y(k) e^{-j \omega T_s k}\right]\left[ \sum_{k=0}^{N-1} u(k) e^{-j \omega T_s k}\right]^{-1}
\end{equation}
where $N$ is the number of data points for each experiment, $u(k) \in \mathbb{R}^{m \times m}$ includes the inputs at instant $k$, $y(k) \in \mathbb{R}^{n \times m}$  the outputs at instant $k$ and $T_s$ is the sampling period. Note that at least $m$ different experiments are needed to extract $G$ from the data (each column of $u(k)$ and $y(k)$ represents respectively the input and the output data from one experiment). We assume that $G(e^{j\omega})$ is bounded in all frequencies except for a set $B_g$ including a finite number of frequencies that correspond to the poles of $G$ on the unit circle. Since the frequency function $G(e^{j\omega})$ is periodic, we consider:
\begin{equation}
	\omega \in \Omega_g=\left \{ \omega \left| -\frac{\pi}{T_s} \leq \omega \leq \frac{\pi}{T_s} \right. \right \} \backslash B_g
\end{equation}

A fixed-structure matrix transfer function controller is considered. The controller is defined as
$K=XY^{-1}$, where $X$ and $Y$ are polynomial matrices in $s$ for continuous-time or in $z$ for discrete-time controller design. This controller structure, therefore, can be used for both continuous-time or discrete-time controllers. The matrix $X$ has the following structure:
\begin{equation}
X=\left[\begin{array}{ccc}X_{11} & \ldots & X_{1n} \\\vdots & \ddots & \vdots \\X_{m1} & \ldots & X_{mn}\end{array}\right] \circ F_x
\end{equation} 
where $X$ and $F_x$ are $m \times n$  polynomial matrices and $\circ$ denotes the element by element multiplication of matrices. The matrix $F_x$ represents the fixed known terms in the controller that are designed to have specific performance, e.g. based on the internal model principle. For discrete-time controllers, we have:
\begin{equation} \label{Xz}
X(z)=X_p z^p+\cdots+X_1 z+X_0
\end{equation}
where $X_i \in \mathbb{R}^{m \times n}$ for $i=0,\ldots,p$ contain the controller parameters. In the same way the matrix polynomial $Y$ can be defined as:
\begin{equation}
Y=\left[\begin{array}{ccc}Y_{11} & \ldots & Y_{1n} \\ \vdots & \ddots & \vdots \\Y_{n1} & \ldots & Y_{nn}\end{array}\right] \circ F_y
\end{equation} 
where $Y$ and $F_y$ are $n \times n$ polynomial matrices. The matrix $F_y$ represents the fixed terms of the controller, e.g. integrators or the denominator of other disturbance models. The set of frequencies of all roots of the determinant of $F_y$ on the stability boundary (imaginary axis for continuous-time controllers or the unit circle for the discrete-time case) is denoted by $B_y$.  

The matrix $Y$ for discrete-time case can be written as:
\begin{equation}  \label{Yz}
Y(z)=I z^p+\cdots+Y_1 z+ Y_0
\end{equation}
where $Y_i \in \mathbb{R}^{n \times n}$ for $i=0,\ldots, p-1$ contain the controller parameters. In order to obtain low-order controllers, a diagonal structure can be considered for $Y$ that makes its inversion and implementation easier too. Note that $Y(e^{j\omega})$ should be invertible for all $\omega \in \Omega=\Omega_g\backslash B_y$. 

The control structure defined in this section is very general and covers centralized, decentralized and distributed control structures. The well-known PID control structure for MIMO systems is also a special case of this structure.

\section{Control Performance}\label{sec3}
%
%
It is shown in this section that classical control performance constraints can be transformed to constraints on the spectral norm of the system and in general can be reformulated as:
\begin{equation}
F^*F-P^*P < \gamma I
\end{equation}
where $F \in \mathbb{C}^{n \times n}$ and $P \in \mathbb{C}^{n \times n}$ are linear in the optimization variables and $(\cdot)^*$ denotes the complex conjugate transpose. This type of constraint is called convex-concave constraint and can be convexified using the 
Taylor expansion of $P^*P$ around $P_c \in \mathbb{C}^{n \times n}$ which is an arbitrary known matrix \citep{DGMD12}:
\begin{equation}
P^*P \approx P_c^*P_c+ (P-P_c)^*P_c+P_c^*(P-P_c)
\end{equation}
It is easy to show that the left hand side term is always greater than or equal to the right hand side term, i.e. :
\begin{equation} \label{Linapprox}
P^*P \geq P^*P_c+P_c^*P- P_c^*P_c
\end{equation}
This can be obtained easily by development of the inequality
$(P-P_c)^*(P-P_c) \geq 0$.

\subsection{$H_\infty$ performance}
Constraints on the infinity-norm of any weighted sensitivity function can be considered. For example, consider the mixed sensitivity problem:
\begin{equation} \label{mixedsens}
\min_K \left\|\begin{array}{c}W_1 S \\W_2 K S\end{array}\right\|_{\infty}
\end{equation}
where $S=(I+GK)^{-1}$ is the sensitivity function, $W_1$ is the performance weight and $W_2$ is the input weight. This problem can be converted to an optimization problem on the spectral norm as:
%
\begin{equation}  \label{spectralcons}
\begin{split}
&\hspace{2cm} \min_K \gamma  \\
&\text{subject to:} \\
&\left[\begin{array}{c}W_1 S \\W_2 K S\end{array}\right]^*
\left[\begin{array}{c}W_1 S\\W_2 K S\end{array}\right] < \gamma I, \qquad \forall \omega \in \Omega
\end{split}
\end{equation}
Note that the argument $e^{j\omega}$ has been omitted for $W_1(e^{j\omega})$, $S(e^{j\omega}), K(e^{j\omega})$ and $W_2(e^{j\omega})$ in order to simplify the notation. The above constraint can be rewritten as:
\begin{multline}
[W_1(I+GK)^{-1}]^*[W_1(I+GK)^{-1}]+\\ [W_2K(I+GK)^{-1}]^*[W_2K(I+GK)^{-1}] < \gamma I
\end{multline}
and converted to a convex-concave constraint as follows:
\begin{multline}
Y^*W_1^* \gamma^{-1} W_1 Y + X^*W_2^* \gamma^{-1} W_2 X \\
- (Y+GX)^*(Y+GX) < 0
\end{multline}
If we denote $P=Y+GX$, using (\ref{Linapprox}), a convex approximation of the constraint can be obtained around $P_c=Y_c+GX_c$ as:
\begin{multline}
Y^*W_1^* \gamma^{-1} W_1 Y + X^*W_2^* \gamma^{-1} W_2 X \\
- P^*P_c-P_c^*P+ P_c^*P_c < 0
\end{multline}
Therefore,  using the Schur complement lemma, the $H_\infty$ mixed sensitivity problem can be represented as the following convex optimization problem with linear matrix inequalities (LMIs):
\begin{equation} \label{MixedHinfty}
\begin{split}
&\hspace{3cm} \min_{X, Y} \gamma \\
&\text{subject to:} \\
& \left[\begin{array}{ccc}P^*P_c+P_c^*P- P_c^*P_c & (W_1Y)^* & (W_2X)^* \\W_1Y & \gamma I & 0 \\W_2X & 0 & \gamma I\end{array}\right] > 0
\end{split}
\end{equation}
for all $\omega \in \Omega$. This convex constraint is a sufficient condition for the spectral constraint in (\ref{spectralcons}) for any choice of an initial controller $K_c=X_cY_c^{-1}$. 

\subsection{$H_2$ performance}
In this section, we show how the $H_2$ control performance can be formulated as a convex optimization problem. We consider the following $H_2$ control performance:
\begin{equation}
	\min_K \|W_1 S\|_2^2
\end{equation}
For a stable closed-loop system, this is equivalent to:
\begin{equation}
\begin{split}
&\hspace{1.5cm}\min_K \int_{-\frac{\pi}{T_s}}^{\frac{\pi}{T_s}} \mbox{trace}[\Gamma(\omega)]d\omega \\
&\text{subject to:}\\
&W_1[(I+GK)^*(I+GK)]^{-1}W_1^* < \Gamma(\omega) \quad \forall \omega \in \Omega
\end{split}
\end{equation}
where $\Gamma(\omega) > 0$ is an unknown matrix function  $\in \mathbb{R}^{n \times n}$. Replacing $K$ with $XY^{-1}$, we obtain:
$$
W_1Y[(Y+GX)^*(Y+GX)]^{-1}Y^*W_1^* < \Gamma(\omega) \quad \forall \omega \in \Omega
$$
which is equivalent to the following matrix inequality:
\begin{equation}
\left[\begin{array}{cc}\Gamma(\omega) & W_1Y \\Y^*W_1^* \, \,& (Y+GX)^*(Y+GX)\end{array}\right] > 0, \quad \forall \omega \in \Omega
\end{equation}
The quadratic part can be linearized using (\ref{Linapprox}) to obtain a linear matrix inequality as:
\begin{equation} \label{H2Prob}
\left[\begin{array}{cc}\Gamma(\omega) & W_1Y \\Y^*W_1^* & \quad P^*P_c+P_c^*P- P_c^*P_c\end{array}\right] > 0, \quad \forall \omega \in \Omega
\end{equation}
{\bf Remark:}
	The unknown function $\Gamma(\omega)$ can be approximated by a polynomial function of finite order as:
	\begin{equation}
	\Gamma(\omega)=\Gamma_0+\Gamma_1 \omega+\cdots+\Gamma_h \omega^h
	\end{equation} 
	In case the constraints are evaluated for a finite set of frequencies $\Omega_N=\{\omega_1,\ldots,\omega_N \}$, $\Gamma(\omega)$ can be replaced with a matrix variable $\Gamma_k$ at each frequency $\omega_k$.

\subsection{Loop shaping}
Assume that a desired loop transfer function $L_d$ is available and that the objective is to design a controller $K$ such that the loop transfer function $L=GK$ is close to $L_d$ in the 2- or $\infty$-norm sense. The objective function for the $\infty$-norm case is to minimize $\|L-L_d\|_{\infty}$ and can be expressed as follows:
\begin{equation}
\begin{split}
& \hspace{2cm}\min_K \gamma \\
&\text{subject to:}\\
&(GK-L_d)^*(GK-L_d) < \gamma I  \qquad \forall \omega \in \Omega 
\end{split}
\end{equation}
Replacing $K$ with $XY^{-1}$ in the constraint, we obtain:
\begin{equation}
(GX-L_dY)^* \gamma^{-1} (GX-L_d Y) - Y^*Y < 0
\end{equation}
Again $Y^*Y$ can be linearized around $Y_c$ using the linear approximation in (\ref{Linapprox}). Thus, the following convex formulation is obtained:
\begin{equation} \label{LSHinfty}
\begin{split}
& \hspace{3cm} \min_{X,Y} \gamma\\
&\text{subject to:}\\
&\left[\begin{array}{cc}Y^*Y_c+Y_c^*Y-Y_c^*Y_c & \quad(GX-L_dY)^* \\GX-L_dY & \gamma I \end{array}\right] > 0 \quad 
\end{split}
\end{equation}
for all $\omega \in \Omega$. 
In a similar way, for minimizing $\| L-L_d  \|_2^2$ the following convex optimization problem can be solved:
\begin{equation} \label{LSH2}
\begin{split}
&\hspace{2cm}\min_{X,Y} \int_{-\frac{\pi}{T_s}}^{\frac{\pi}{T_s}} \mbox{trace}[\Gamma(\omega)]d\omega\\
&\text{subject to:}\\
&\left[\begin{array}{cc}Y^*Y_c+Y_c^*Y-Y_c^*Y_c & \quad(GX-L_dY)^* \\GX-L_dY & \Gamma(\omega)  \end{array}\right] > 0 \quad 
\end{split}
\end{equation}
for all $\omega \in \Omega$. Note that the resulting loop shaping controller does not necessarily guarantee the closed-loop stability. This will be discussed in the next section, where the stability conditions will be developed.


\section{Robust Controller Design}\label{sec4}

\subsection{Stability analysis}\label{sec4.1}

The stability of the closed-loop system is not necessarily guaranteed even if the spectral norm of a weighted sensitivity function is bounded.  In fact, an unstable system with no pole on the stability boundary has a bounded spectral norm. In this section, we show that the closed-loop stability can be guaranteed if some conditions in the linearization of the constraints are met. More precisely, the initial controller $K_c=X_cY_c^{-1}$ plays an important role in guaranteeing the stability of the closed-loop system with the resulting controller $K$. Our stability analysis is based on the generalized Nyquist stability criterion for MIMO systems that is recalled here for discrete-time systems. Note that the results are also straightforwardly applicable to the continuous-time case by modifying the Nyquist contour.
\begin{theorem} {\bf(Nyquist stability theorem)}
	The closed-loop system with the plant model $G(z)$ and the controller $K(z)$ is stable if and only if the Nyquist plot of $\det(I+G(z)K(z))$
\begin{enumerate}
	\item makes $N_G+N_K$ counterclockwise encirclements of the origin, where $N_G$ and $N_K$ are, respectively, the number of poles of $G(z)$ and $K(z)$ on the exterior of the unit circle, and
	\item does not pass through the origin.
\end{enumerate}
\end{theorem}
The Nyquist plot is the image of $\det(I+GK)$ as $z$ traverses the Nyquist contour (the unit circle) counterclockwise. We assume that the Nyquist contour has some small detours around the poles of $G(z)$ and $K(z)$ on the unit circle.
 
\begin{definition}  Let $\wno\{F(z)\}$ be the winding number, in the counterclockwise sense, of the image of $F(z)$ around the origin when $z$ traverses the Nyquist contour with some small detours around the poles of $F(z)$ on the unit circle.  
\end{definition}
Since the winding number is related to the phase of the complex function, we have the following properties:
\begin{align}
&\wno \{F_1(z)F_2(z)\}=\wno \{F_1(z)\}+\wno\{F_2(z)\} \\
&\wno \{F(z)\}=- \wno \{F^*(z)\} \\
&\wno \{F(z)\}=- \wno \{F^{-1}(z)\} 
\end{align}

\begin{theorem} \label{theorem:stability}
Given a plant model $G$, an initial stabilizing controller $K_c=X_cY_c^{-1}$ with $\det(Y_c)\neq 0, \forall \omega \in \Omega$, and feasible solutions $X$ and $Y$ to the following LMI, 
\begin{equation}  \label{StabCon}
	(Y+GX)^*(Y_c+GX_c)+(Y_c+GX_c)^*(Y+GX) > 0 
\end{equation}
for all $\omega \in \Omega$, then the controller $K=XY^{-1}$ stabilizes the closed-loop system if 
\begin{enumerate}
	\item $\det(Y) \neq 0, \forall \omega \in \Omega$.
	\item The initial controller $K_c$ and the final controller $K$ share the same poles on the stability boundary, i.e. $\det(Y)=\det(Y_c)=0, \forall \omega \in B_y$.
	\item The order of $\det(Y)$ is equal to the order of $\det(Y_c)$.
\end{enumerate}
\end{theorem}

{\it Remark:} Note that the condition in (\ref{StabCon}) is always met when a convexified $H_\infty$ or $H_2$ control problem has a feasible solution because we have $P^*P_c+P_c^*P >0$ in (\ref{MixedHinfty}) and (\ref{H2Prob}).

\textbf{Proof:} The proof is based on the Nyquist stability criterion and the properties of the winding number. The winding number of the determinant of $P^*(z)P_c(z)$ is given by: 
\begin{align} \label{Eq:wno}
	\wno\{\det(P^*P_c) \} =& \, \wno\{\det(P^*)\}+\wno\{\det(P_c)\}  \nonumber \\
	=&-\wno\{\det(I+GK)\det(Y)\} \nonumber  \\ & +\wno\{\det(I+GK_c)\det(Y_c)\} \nonumber \\
	=&-\wno\{\det(I+GK)\} \nonumber \\& -\wno\{\det(Y)\}+\wno\{\det(Y_c)\} \nonumber \\&+\wno\{\det(I+GK_c)\}
\end{align}
Note that the phase variation of $\det(P^*P_c)$ for the small detour in the Nyquist contour is zero, if Condition 2 of the theorem is satisfied. In fact for each small detour, the Nyquist plot of  $\det(I+GK)$ and $\det(I+GK_c)$ will have the same phase variation because $K$ and $K_c$ share the same poles on the unit circle. As a result,
the winding number of $\det(P^*P_c)$ can be evaluated on $\Omega$ instead of the Nyquist contour. On the other hand, the condition in (\ref{StabCon}) implies that $P^*(e^{j\omega})P_c(e^{j\omega})$ is a non-Hermitian positive definite matrix in the sense that :
\begin{equation}
	\Re \{x^*P^*(e^{j\omega})P_c(e^{j\omega})x\} > 0  \qquad  \forall x\neq 0 \in \mathbb{C}^n 
\end{equation}
and $\forall \omega \in \Omega$. This, in turn, means that all eigenvalues of  $P^*(e^{j\omega})P_c(e^{j\omega})$, denoted $\lambda_i(\omega)$ for $i=1,\ldots,n$, have positive real parts at all frequencies \citep{ZLHL10}: 
\begin{equation} \label{eq:eigposreal}
	\Re\{\lambda_i (\omega) \} > 0 \qquad \forall \omega \in \Omega,  i=1,\ldots,n
\end{equation}
Therefore, $\lambda_i(\omega)$ will not pass through the origin and not encircle it (i.e. its winding number is zero). As a result, since the determinant of a matrix is the product of its eigenvalues, we have:
$$\wno\{\det(P^*P_c)\} = \wno \left\{ \prod_{i=1}^n \lambda_i \right\} =\sum_{i=1}^n  \wno \{\lambda_i\}=0$$
 Since $K_c$ is a stabilizing controller, based on the Nyquist theorem
$\wno\{\det(I+GK_c)\}=N_G+N_{K_c}$. Furthermore, according to the argument principle $\wno\{ \det(Y) \}=\delta-N_K$ and $\wno\{ \det(Y_c) \}=\delta-N_{K_c}$, where $\delta$ is the order of $\det(Y)$ and $\det(Y_c)$ according to Condition 3. Now using (\ref{Eq:wno}), we obtain:
\begin{align}
\wno\{\det(I+GK)\}= & \wno\{\det(I+GK_c)\} \nonumber \\ &-\wno\{\det(Y)\}+\wno\{\det(Y_c)\} \nonumber \\
=&N_G+N_K
\end{align}
which shows that Condition 1 of the Nyquist theorem is met. 
We can see from (\ref{eq:eigposreal}) that 
\begin{equation}
\det(P^*P_c) = \prod_{i=1}^n \lambda_i(\omega) \neq 0  \qquad \forall \omega \in \Omega
\end{equation}
Therefore, $\det(P)=\det(I+GK) \det(Y) \neq 0$ and the Nyquist plot of $\det(I+GK)$ does not pass through the origin and Condition 2 of the Nyquist theorem is also satisfied. 
 {\hfill \ensuremath{\blacksquare}}

{\bf Remark 1:} A necessary and sufficient condition for $\det(Y)\neq0$ is $Y^*Y > 0$.
Since this constraint is concave, it can be linearized to obtain the following sufficient LMI:
\begin{equation}
Y^*Y_c+Y_c^*Y-Y_c^*Y_c > 0
\end{equation}
This constraint can be added to the optimization problem in (\ref{MixedHinfty}) in order to guarantee the closed-loop stability for the mixed sensitivity problem. For the loop-shaping problems in (\ref{LSHinfty}) and in (\ref{LSH2}), this condition is already included in the formulation. Therefore, for guaranteeing the closed-loop stability, the condition in (\ref{StabCon}) should be added. This condition can be added directly or by considering an additional $H_2$ or $H_\infty$ constraint on a closed-loop sensitivity function.

{\bf Remark 2:} In practice, condition 3 of Theorem~\ref{theorem:stability} is not restrictive. Any initial controller of lower order than the final controller can be augmented by adding an appropriate number of zeros and poles at the origin in $X_c$ and $Y_c$, thus satisfying the condition without affecting the initial controller.

\subsection{Multimodel uncertainty}

The case of robust control design with multimodel uncertainty is very easy to incorporate in the given framework.
Systems that have different frequency responses in $q$ different operating points can be represented by a multimodel uncertainty set:
\begin{equation}  \label{MMuncertainty}
	\mathcal{G}(e^{j\omega})=\{ G_1(e^{j\omega}), {G}_2(e^{j\omega}), \ldots, {G}_q(e^{j\omega}) \}
\end{equation}
Note that the models may have different orders and may contain the pure input/output time delay.

This can be implemented by formulating a different set of constraints for each of the models. Let $P_i = Y + G_i X$ and $P_{c_i} = X_c + G_i Y_c$. Again taking the mixed sensitivity problem as an example, the formulation of this problem including the stability constraint would be: 
\begin{align} 
&\hspace{3cm} \min_{X,Y} \gamma  \nonumber\\
&\text{subject to:}  \nonumber\\
& \left[\begin{array}{ccc}P_i^*P_{c_i}+P_{c_i}^*P_i- P_{c_i}^*P_{c_i} & (W_1Y)^* & (W_2X)^* \\W_1Y & \gamma I & 0 \\W_2X & 0 & \gamma I \end{array}\right] > 0 \nonumber\\
& Y^*Y_c+Y_c^*Y-Y_c^*Y_c  > 0 \\
& \text{for } i=1,\ldots,q \quad ; \quad \forall \omega \in \Omega \nonumber
\end{align}

\subsection{Frequency-domain uncertainty}

The frequency function may be affected by the measurement noise. In this case, the model uncertainty can be represented as :
\begin{equation}  \label{additive}
	\tilde{G}(e^{j\omega})=G(e^{j\omega})+W_1(e^{j\omega}) \Delta W_2(e^{j\omega})
\end{equation}
where $ \Delta$ is the unit ball of matrices of appropriate dimension and $W_1(e^{j\omega})$ and $W_2(e^{j\omega})$ are known complex matrices that specify the magnitude of and directional information about the measurement noise. A convex optimization approach is proposed in \cite{HSB02} to compute the optimal uncertainty filters from the frequency-domain data. The system identification toolbox of Matlab provides the variance of $G_{ij}(e^{j\omega})$ (the frequency function between the $i$-th output and the $j$-th input) from the estimates of the noise variance that can be used for computing $W_1$ and $W_2$.  

The robust stability condition for this type of uncertainty is \citep{Zho98}:
$\| W_2 K S W_1 \|_\infty < 1$.
If we assume that $W_1(e^{j\omega})$ is invertible for all $\omega \in \Omega$ (i.e. it has no pole on the unit circle), then a set of robustly stabilizing controllers can be given by the following spectral constraints:
 \begin{align} 
& \left[\begin{array}{ccc}P^*P_{c}+P_{c}^*P- P_{c}^*P_{c} & (W_2X)^* \\W_2X & I \end{array}\right] > 0 \\
& Y^*Y_c+Y_c^*Y-Y_c^*Y_c  > 0 \quad ; \quad  \forall \omega \in \Omega \nonumber
\end{align}
where $P=W_1^{-1}(Y+GX)$ and $P_c=W_1^{-1}(Y_c+GX_c)$.

\section{Implementation Issues} \label{sec5}

%

\subsection{Frequency gridding}

The optimization problems formulated in this paper contain an infinite number of constraints (i.e. $\forall \omega \in \Omega$) and are called semi-infinite problems. A common approach to handle this type of constraints is to choose a reasonably large set of frequency samples $\Omega_N = \left\{\omega_1, \ldots, \omega_N \right\}$ and replace the constraints with a finite set of constraints at each of the given frequencies. As the complexity of the problem scales linearly with the number of constraints, $N$ can be chosen relatively large without severely impacting the solver time. The frequency range $\left[ 0, \pi/T_s \right]$ is usually gridded  logarithmically-spaced. Since all constraints are applied to Hermitian matrices, the constraints for the negative frequencies between $-\pi/T_s$ and zero will be automatically satisfied. In some applications with low-damped resonance frequencies, the density of the frequency points can be increased around the resonant frequencies. An alternative is to use a randomized approach for the choice of the frequencies at which the constraints are evaluated \citep{ATL10}. 

Taking the mixed sensitivity problem as an example, the sampled problem would be: 
\begin{equation} \label{mixedgrides}
\begin{split}
&\hspace{3cm} \min_{X,Y} \gamma \\
&\text{subject to:} \\
& \left[\begin{array}{ccc}P^*P_c+P_c^*P- P_c^*P_c & (W_1Y)^* & (W_2X)^* \\W_1Y & \gamma I & 0 \\W_2X & 0 & \gamma I \end{array}\right](e^{j \omega}) > 0 \\[5pt]
& \left[\begin{array}{c} Y^*Y_c+Y_c^*Y-Y_c^*Y_c \end{array}  \right] (e^{j \omega}) > 0  \quad ; \quad  \omega \in \Omega_N
\end{split}
\end{equation}

\subsection{Initial controller}

The stability condition presented in Theorem~\ref{theorem:stability} requires a stabilizing initial controller $K_c$ with the same poles on the stability boundary (the unit circle) as the desired final controller. For a stable plant, a stabilizing initial controller can always be found by choosing: 
\begin{equation}
	\left[ X_{c,1}, \ldots, X_{c,p} \right] = 0 , \quad X_{c,0} = \epsilon I
\end{equation}
with $\epsilon$ being a sufficiently small number. Furthermore, the parameters of $Y_c$
should be chosen such that $\det(Y_c) \neq 0$ for all $\omega \in \Omega$. This can be achieved by choosing $Y_c$ such that all roots of $\det(Y_c)=0$ lie at zero, with $F_y$ containing all the poles on the unit circle of the desired final controller. For example, to design a controller with integral action in all outputs, $Y_c=z^p(z-1)I$ can be considered. Alternatively, if a working controller has already been implemented, it can be used as the initial controller. 

When choosing an initial controller whose performance is far from the desired specifications, it may occur that either the optimization problem has no feasible solution, or that the solver runs into numerical problems which lead to an infeasible solution. These problems can often be resolved by two approaches:
\begin{description}
	\item[Re-initialization:] The initial controller can be changed with a systematic approach for stable plants by solving the following optimization problem using a nonlinear optimization solver with random initialization: 
\begin{equation} 
\begin{split}
&\hspace{2cm} \max_{X,Y} \, a \\
& \text{subject to:} \\
& \Re \left \{\det(I + GXY^{-1}) \right \} \geq a \qquad \forall \omega \in \Omega_N 
\end{split}
\end{equation}
Any solution to the above optimization problem will be a stabilizing controller if the optimal value of $a$ is greater than -1. The problem can be solved multiple times with different random initialization to generate a set of initial stabilizing controllers, which can be used to initialize the algorithm. 

	\item[Relaxation:] We can relax or even remove some of the constraints. The relaxed optimization problem is then solved and the optimal controller is used to initialize the non-relaxed problem. As this new controller is comparatively close to the final solution, the issue is often solved with this approach.
\end{description}

Since this work focuses on data-driven control design, for unstable plants it is reasonable to assume that a stabilizing controller has been available for data acquisition, and can thus be used as the initial controller.

It should be mentioned that the design of fixed-structure controllers in a model-based setting also requires an initialization with a stabilizing controller, which is usually integrated in the workflow. The methods based on non-smooth optimization like {\it hinfstruct} in Matlab \citep{AN06} or the public-domain toolbox HIFOO \citep{BHLO06} use a set of randomly chosen stabilizing controllers for initialization and take the best result. This set is constructed by solving a non-convex optimization problem that minimizes the maximum eigenvalue of a closed-loop transfer function. Other model-based approaches use an initial stabilizing controller to convert the bilinear matrix inequalities to LMIs and solve it with convex optimization algorithms. Therefore, from this point of view, our data-driven approach is subject to the same restrictions as the state-of-the-art approaches for fixed-structure controller design in a model-based setting.

\subsection{Iterative algorithm}

Once a stabilizing initial controller is found, it is used to formulate the optimization problem. Any LMI solver can be used to solve the optimization problem and calculate a suboptimal controller $K$ around the initial controller $K_c$. As we are only solving an inner convex approximation of the original optimization problem, $K$ depends heavily on the initial controller $K_c$ and the performance criterion can be quite far from the optimal value. The solution is to use an iterative approach that solves the optimization problem multiple times, using the final controller $K$ of the previous step as the new initial controller $K_c$. This choice always guarantees closed-loop stability (assuming the initial choice of $K_c$ is stabilizing). Since the objective function is non-negative and non-increasing, the iteration converges to a local optimal solution of the original non-convex problem \citep{YR03}. The iterative process can be stopped once the change in the performance criterion is sufficiently small.

\section{Simulation Results}\label{sec6}

As an example, the mixed sensitivity problem for low-order continuous-time controllers is considered. 10 plants are drawn from the Compleib library \citep{Lei06}. For comparison, the achieved performance is compared with the results obtained using {\it hinfstruct} and HIFOO. Parametric plant models are used in this example in order to enable comparison with state-of-the-art methods. However, it should be noted that, as our method is data-driven, only the frequency responses of the plants are required for the controller design.

The objective is to solve the mixed sensitivity problem by minimizing the infinity-norm of (\ref{mixedsens}), where $W_2 = I$ and $W_1 = (a_k s + 10)/(a_k s + 1)$
with $a_k$ being chosen based on the bandwidth of the plant. Then, the optimization problem in (\ref{mixedgrides}) is formed with $N=100$ logarithmically spaced frequency points in the interval $\left[0.01, 500 \right]$ rad/s, where $500$ is much larger than the bandwidth of all plants. 
A second-order controller $K(s) = X(s)Y(s)^{-1}$ is chosen as follows: 
$$
	X(s) = X_2 s^2 + X_1 S + X_0 \quad , \quad
	Y(s) = I s^2 + Y_1 s + Y_0
$$
where $Y_i$ is a diagonal matrix in order to obtain a low-order controller. To have a fair comparison, the same method as in HIFOO is used to find a stabilizing initial controller. The method uses a non-convex approach to minimize the maximum of the spectral abscissa of the closed-loop plant, and yields a stabilizing static output feedback controller $K_{\text{SOF}}$. In order to satisfy Condition~3 of Theorem~2, the order of $Y_c$ is increased without changing the initial controller :
\begin{equation}
	X_c(s) = (s+1)^2 K_{\text{SOF}} \ , \ Y_c(s) = (s+1)^2 I
\end{equation}
The names of the chosen plants in Compleib, the design parameters and the obtained norms are shown in Table~\ref{tab:exparams}. For comparison, the mixed sensitivity problems are also solved for a second-order state-space controller using HIFOO and {\it hinfstruct} with 10 random starts. It can be seen that the data-driven method generally achieves about the same or a lower norm. The superior results can be attributed to the fact that the controller structure is of matrix polynomial form, which has more parameters than a state-space controller of the same order. 

The solver time of one iteration step depends almost linearly on the number of points used for the frequency gridding. It is also interesting to note that the controller order has a minimal impact on the solver time, making the algorithm well-suited for the design of higher-order controllers. The number of iterations until convergence mostly depends on the choice of the initial controller and a solution is generally reached in less than 25 iterations.

\begin{table}
	\centering	
	\caption{Comparison of optimal mixed sensitivity norms for 10 plants from Compleib}
	\label{tab:exparams}
    	\begin{tabular}{c|c|c|c|c}
		\hline 
    		Plant Name & $a_k$ & data-driven & {\it hinfstruct} & HIFOO \\
		\hline 
				AC1 & 	10 &		1.90 &	2.30 &	2.38 \\
				HE1 & 	1 & 	1.37 & 	1.36 &  1.36 \\
				HE2 & 	10 & 	3.08 & 	3.36 & 	3.55 \\
				REA2 & 	1 & 	3.00 & 	2.96& 	2.96 \\
				DIS1 & 	1 & 	7.27 & 	7.31 & 	7.34 \\
				TG1 & 	0.1 &	9.54 & 	8.89& 	9.75 \\
				AGS & 	1 & 	2.14 & 	2.16 & 	2.16 \\
				BDT2 & 	1 & 	9.93 & 9.93& 	9.94 \\
				MFP & 	1 & 	6.08 & 	7.23 & 	7.17 \\
				IH & 		1 & 	4.83 & 	10.01 & 28.73 \\   		
    		\hline 
    	\end{tabular} 
\end{table}

\section{Further Simulation Results}

In this section, three additional examples are presented to demonstrate the applicability of the method. Note that in the first two examples, for the sake of comparison with model-based methods, a parametric model of the plant is given. However, this parametric model is not used in the controller design and only its frequency response $G(j\omega)$ is employed. For each example, the optimization problem was formulated in Matlab using Yalmip~\cite{yalmip}, and solved with Mosek~\cite{Mosek}. 

\subsection{Fixed-structure controller design}

The first example is drawn from Matlab's Robust Control Toolbox and treats the control design for a 9th-order model of a head-disk assembly in a hard-disk drive. In the Matlab example, {\it hinfstruct} is used to design a robust controller such that a desired open-loop response is achieved while satisfying a certain performance measure. We will show that an equivalent controller of the same order can be designed using the method presented in this paper.

%
\begin{figure}[t]
    \centering
        \includegraphics[width=0.9\columnwidth]{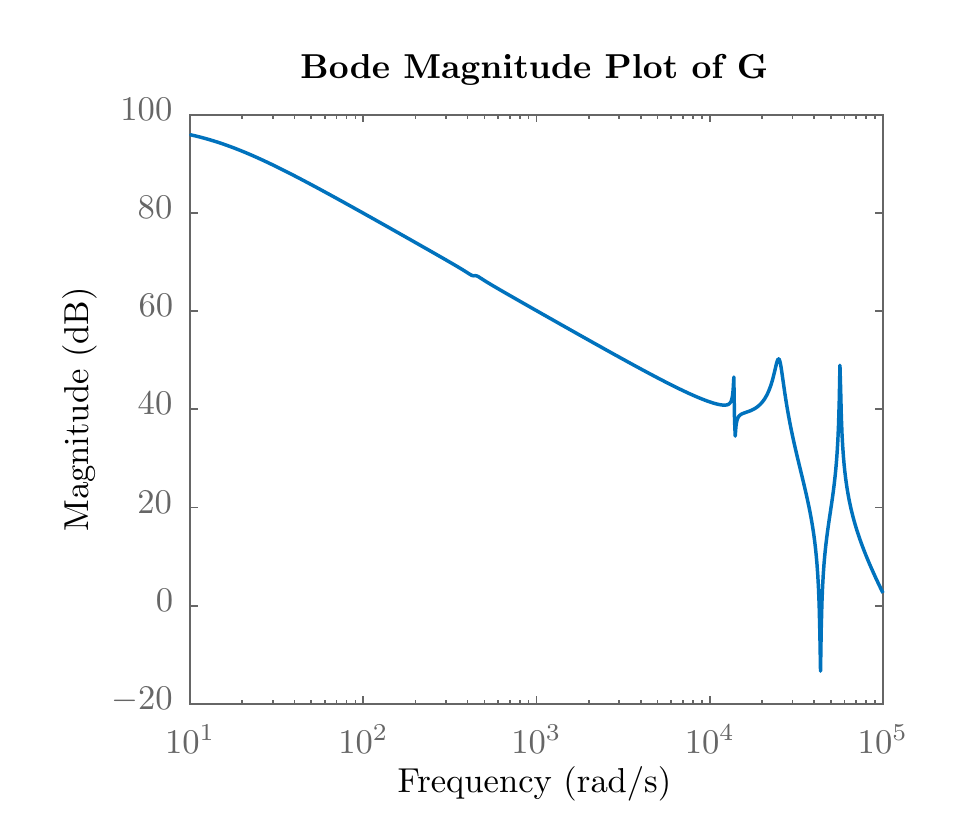}
        \caption{Bode magnitude plot of the plant used in example 1.}
				\label{fig:hda-bode}
\end{figure}

\begin{figure}
	\centering
		\includegraphics[width=0.9\columnwidth]{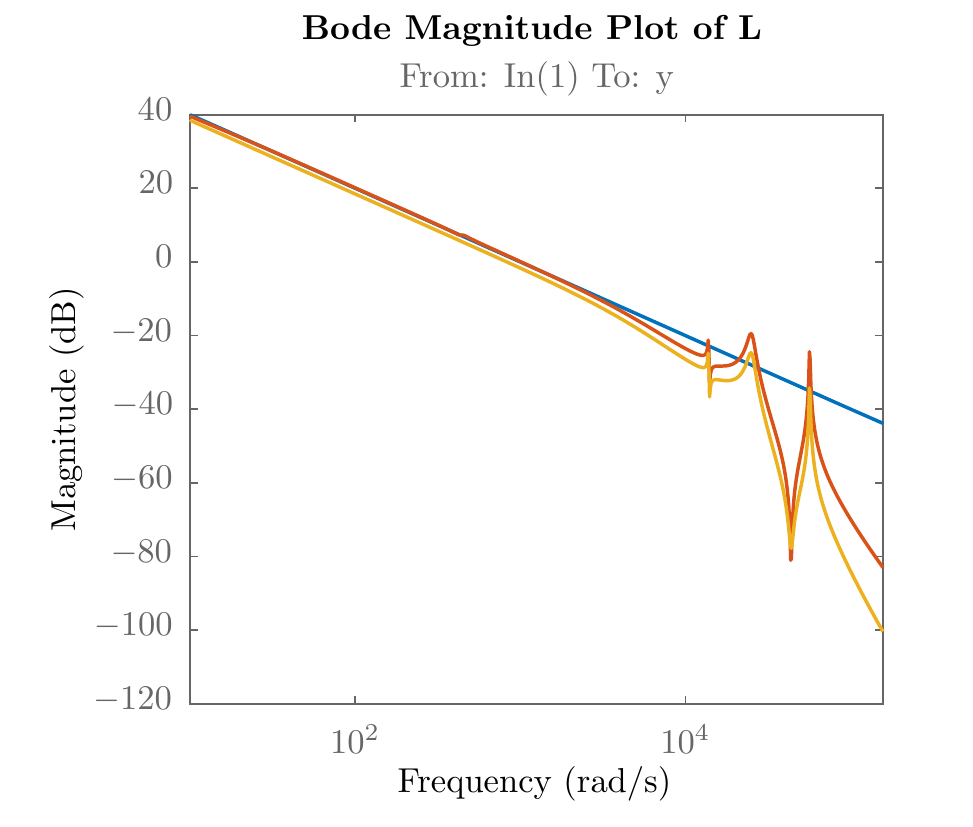}
	\caption{Comparison of the open-loop transfer functions; The blue line: desired open-loop, the red line: proposed method, the yellow line: {\it hinfstruct} controller.}
	\label{fig:ex1_bodemag_L_comparison}
\end{figure}

The bode magnitude plot of the plant is shown in Fig.~\ref{fig:hda-bode}. The desired open-loop transfer function is given by: 
\begin{equation}
	L_d(s) = \frac{s + 10^6}{1000s + 1000} 
\end{equation}
Additionally, a constraint on the closed-loop transfer function is introduced to increase the robustness and performance: 
$\left\| W_1 T \right\|_{\infty} \leq 1$ and $W_1 = 1$.
To stay in line with the data-driven focus of this paper, we choose to design a discrete-time controller with the same order as the continuous-time controller given in the Matlab example: 
\begin{equation}
	K(z) = \frac{X_2 z^2 + X_1 z + X_0}{(z-1) (z + Y_0)}
\end{equation}
Since the plant is stable, an initial controller is easily found by setting $X_1,X_2,Y_0$ to zero and choosing a small enough value for $X_0$. This results in the following initial controller: 
\begin{equation}
	K_{c}(z) = \frac{ 10^{-6}}{z^2 - z}
\end{equation}
Note how the pole on the unit circle introduced by the integrator is also included in the initial controller. Then the problem is formulated as an $H_2$ loop shaping problem. The semi-infinite formulation is sampled using 1000 logarithmically spaced frequency points in the interval $\Omega_N = \left[ 10, 5\times10^4 \pi \right]$ (the upper limit being equal to the Nyquist frequency). The semi-definite problem is as follows: 
\begin{align}
&\hspace{2cm}\min \sum_{k = 1}^{N} \mbox{trace}[\Gamma_k]  \nonumber\\
&\text{subject to:} \nonumber\\
&\left[\begin{array}{cc}Y^*Y_c+Y_c^*Y-Y_c^*Y_c & \quad(GX-L_dY)^* \\G X-L_dY & \Gamma_k \end{array}\right] (j \omega_k) > 0 \nonumber \\
&\left[\begin{array}{cc} P^*P_c + P_c^*P - P_c^*P_c & (W_1GX)^* \\ W_1GX & I \end{array}\right] (j \omega_k) > 0 \\
& k=1,\ldots,N \nonumber
\end{align}
The algorithm converges within 10 iterations to a final, stabilizing controller that satisfies the closed-loop constraint and has the following parameters:
\begin{equation}
	K(z) = 10^{-4} \frac{2.287z^2 - 3.15 z + 0.8631}{(z-1)(z - 0.8598)}
\end{equation}
Fig.~\ref{fig:ex1_bodemag_L_comparison} shows a comparison of the desired open-loop transfer function and the results produced by our method as well as the controller calculated in the Matlab example using {\it hinfstruct}. It can be seen that the result is very similar to the result generated by {\it hinfstruct}, with our result being closer to the desired transfer function at lower frequencies. This is especially noticeable when comparing the sum of the trace of $\Gamma_k$ in the objective function, with our solution achieving a value that is around 30 times smaller. 

\subsection{Mixed sensitivity problem}

In this example the mixed sensitivity problem of a 3$\times$3 MIMO continuous-time plant model is considered. The global optimal solution to this problem with a full-order controller can be obtained via Matlab using \textit{mixsyn}. The plant is taken from the first example in \cite{SOW10} and has the following transfer function: 
\begin{equation}
	G(s) = \left[\begin{array}{ccc}
		\frac{1}{s+1} & \frac{0.2}{s+3} & \frac{0.3}{s+0.5} \\
		\frac{0.1}{s+2} & \frac{1}{s+1} & \frac{1}{s+1} \\
		\frac{0.1}{s+0.5} & \frac{0.5}{s+2} & \frac{1}{s+1}
	\end{array} \right]
\end{equation}
The objective is to solve the mixed sensitivity problem by minimizing the infinity-norm of (\ref{mixedsens}), where 
the weighting transfer functions are also taken from \cite{SOW10}:
\begin{equation}
	W_1 = \frac{s+3}{3s + 0.3} I \ , \ W_2 = \frac{10s + 2}{s+40} I
\end{equation}
In this example we design a continuous-time controller to show that the developed frequency-domain LMIs in this paper can be used directly to design continuous-time controllers. The controller transfer function matrix is defined as $K(s) = X(s)Y^{-1}(s)$, where: 
\begin{align}
	X(s) &= X_p s^p + \ldots +  X_1 s + X_0  \\ 
	Y(s) &= I s^p + \ldots + Y_1 s + Y_0  
\end{align}
and $p$ is the controller order. The optimization problem is sampled using $N=1000$ logarithmically spaced frequency points in the interval $\Omega_N = \left[10^{-2},10^2\right]$, resulting in the following optimization problem : 
\begin{align}
&\hspace{2.5cm}\min_{X,Y} \gamma \nonumber \\
&\text{subject to:}  \nonumber\\
&\left[\begin{array}{ccc} P^*P_c + P_c^*P - P_c^*P_c & (W_1Y)^* & (W_2X)^* \\ 
												 W_1Y & \gamma I & 0 \\
												 W_2X & 0 & \gamma I
			 \end{array}\right] (j \omega_k) > 0 \nonumber\\
&\left[ Y^*Y_c+Y_c^*Y-Y_c^*Y_c \right] (j \omega_k) > 0 \\
&k=1,\ldots,N \nonumber
\end{align}
Since the plant is stable, an initial controller is found by setting the poles of the controller to $-1$, i.e. $Y_c=(s+1)^pI$ and choosing $X_0 = I, \left\{X_1,\ldots,X_p\right\} = 0$.

\begin{figure}
	\centering
	\includegraphics[width=0.75\columnwidth]{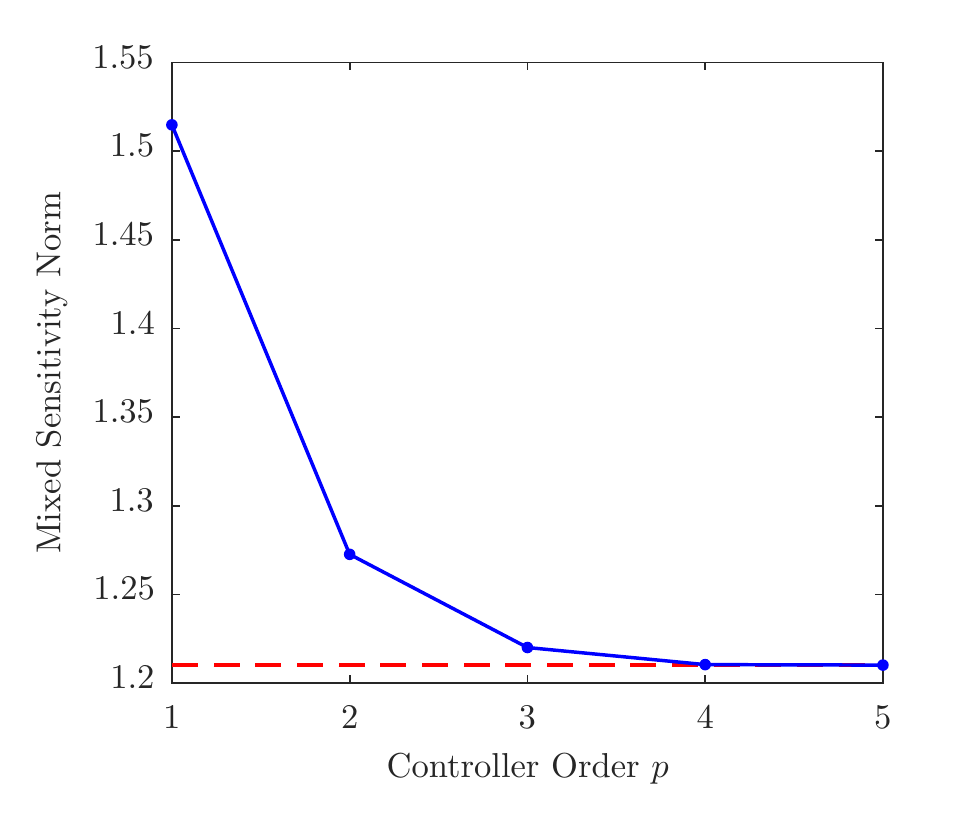}
	\caption{Plot of the mixed sensitivity norm for different controller orders $p$. The dashed red line shows the globally optimal value obtained by \textit{mixsyn}.}
	\label{fig:ex2_gamma}
\end{figure}

The problem is then solved for controller orders $p$ from 1 to 5, with the algorithm converging within 3 to 6 iterations. The value of the obtained norm is shown in Fig.~\ref{fig:ex2_gamma}. The number of design parameters is equal to $(2p+1) \times 9$. The figure also shows the globally optimal norm for a full-order state-space controller with 289 design parameters obtained through {\it mixsyn}. It can be seen that already for $p=3$ a very good value is achieved with the following controller parameters: 
\resizebox{\columnwidth}{!}{
  \begin{minipage}{\linewidth}
\begin{align*}
	X(s) = &\left[ \begin{array}{ccc} \nonumber
		0.0794 & 0.0041 & -0.0032 \\
		0.0091 & 0.1076 & -0.0421 \\
		0.0131 & 0.031 & 0.0986
	\end{array} \right] s^3 + 
	\left[ \begin{array}{ccc}
		4.5304 & -0.6974 & -0.8464 \\
		-0.5345 & 3.2929 & -2.3889 \\
		-0.3737 & -0.1412 & 3.421
	\end{array} \right] s^2 \\ \nonumber
	+ &\left[ \begin{array}{ccc}
		9.0896 & -3.4091 & -2.6272 \\
		2.2293 & 4.0883 & -3.1235 \\
		-3.0827 & -0.3391 & 3.4927
	\end{array} \right] s  + 
  \left[ \begin{array}{ccc}
		2.0218 & -1.0874 & -1.6883 \\
		2.4056 & 1.7292 & -0.6611 \\
		-1.0974 & -0.1376 & 1.8895
	\end{array} \right] \\
Y(s) = &\left[ \begin{array}{ccc} \nonumber
		1 & 0 & 0 \\
		0 & 1 & 0 \\
		0 & 0 & 1
	\end{array} \right] s^3 + 
	\left[ \begin{array}{ccc}
		5.1556 & -1.1562 & -0.5595 \\
		-0.5993 & 1.9965 & -0.6899 \\
		-0.9489 & -0.6155 & 2.2864
	\end{array} \right] s^2 \\ 
	+ &\left[ \begin{array}{ccc}
		2.444 & -1.2479 & -0.7046 \\
		0.729 & 1.427 & 0.0589 \\
		-0.9949 & -0.5552 & 1.1323
	\end{array} \right] s  + 
  \left[ \begin{array}{ccc}
		0.1514 & -0.1487 & -0.1067 \\
		0.2084 & 0.1941 & 0.1491 \\
		-0.0116 & -0.0029 & 0.1791
	\end{array} \right] \\		\nonumber
\end{align*}
\end{minipage}
}
For $p=5$, with only 99 design parameters the global optimum is achieved. This example shows that the proposed method is able to reach the global optimum value of the mixed sensitivity norm for a general MIMO transfer function while having a significantly lower number of design parameters than the classical state-space methods. It also yields good results for lower-order controllers and does not require a parametric model. 

\subsection{Data-driven control of a gyroscope}

For the third example, we design a data-driven, robust multivariable controller with multimodel uncertainty to control the gimbal angles of a gyroscope. We then apply the controller on an experimental setup to validate the results. 

\subsubsection{Experimental setup}

\begin{figure}
	\centering
		\includegraphics[width=0.6\columnwidth]{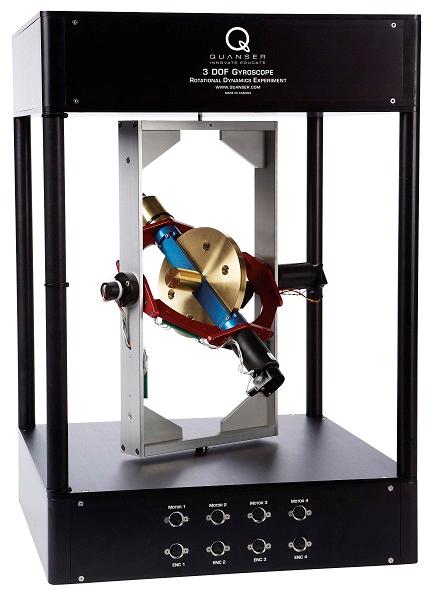}
	\caption{The gyroscope experimental setup by Quanser.}
	\label{fig:Gyroscope}
\end{figure}

The experiment was conducted on a 3 DOF gyroscope setup built by Quanser (see Fig.~\ref{fig:Gyroscope}). The system consists of a disk mounted inside an inner blue gimbal, which is in turn mounted inside an outer red gimbal. The entire structure is supported by the rectangular silver frame. The disk, both gimbals and the frame can be actuated about their respective axis through electric motors, and their angular positions can be measured using high resolution optical encoders. For this experiment, the position of the silver frame is mechanically fixed in place. The control objective is to achieve a good tracking performance on the angular positions of the blue and red gimbal and to minimize the coupling between the axes. The dynamics of the system change depending on the angular velocity of the disk, which is included in the control design as a multimodel uncertainty.

\subsubsection{Plant identification}

\begin{figure}
	\centering

\unitlength .80mm
\linethickness{0.5pt}

\begin{picture}(110.00,35.00)(15.00, -5.00)

\put(33.00,10.00){\framebox(15.00,10.00)[cc]{$K$}}
\put(60.00,10.00){\framebox(15.00,10.00)[cc]{$K_{\text{fl}}$}}
\put(85.00,10.00){\framebox(15.00,10.00)[cc]{$G_m$}}

\put(75.00,15.00){\vector(1,0){10}}
\put(100.00,15.00){\vector(1,0){25}}
\put(105.00,15.00){\line(0,-1){10.00}}
\put(105.00,5.00){\line(-1,0){37.50}}
\put(67.50,5.00){\vector(0,1){5}}

\put(48.00,15.00){\vector(1,0){12}}

\put(25.00,15.00){\circle{4.00}}
\put(115.00,15.00){\line(0,-1){20.00}}
\put(115.00,-5.00){\line(-1,0){90.00}}
\put(25.00,-5.00){\vector(0,1){18.00}}
\put(27.00,15.00){\vector(1,0){6.00}}
\put(15.00,15.00){\vector(1,0){8.00}}

\put(17.00,18.00){\makebox(0,0)[cc]{\small $\theta^*$}}
\put(53.00,18.00){\makebox(0,0)[cc]{\small $\theta_u$}}
\put(120.00,18.00){\makebox(0,0)[cc]{\small $\theta$}}
\put(27.00,12.00){\makebox(0,0)[cc]{-}}

\color{red}
\put(57.00,0.00){\dashbox(53.00,25.00)[cc]}
\put(82.50,28.00){\makebox(0,0)[cc]{$G$}}

\end{picture}
	\caption{Block diagram of the cascaded controller structure of the gyroscope.}
	\label{fig:gyro_controlarchitecture}
\end{figure}
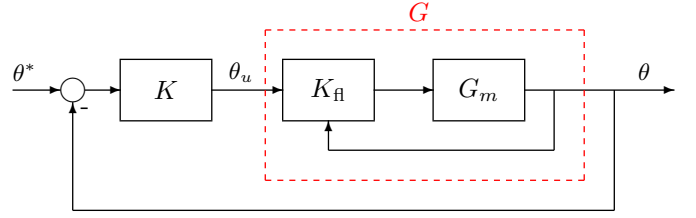

The gyroscope is a strongly nonlinear system, and linear control design methods only achieve good performance in a small range around the operation points. In order to improve this range, a cascaded control architecture was chosen with a feedback linearization forming the inner loop (see Fig.~\ref{fig:gyro_controlarchitecture}). The block $G_m$ is the real plant, $K_{fl}$ is the feedback linearization and $K$ is the controller to be designed. $\theta = \left[ \theta_{b}, \theta_{r}\right], \theta^* = \left[ \theta_{b}^*, \theta_{r}^*\right]$ are vectors containing the measured and desired blue and red gimbal angles. $\theta_u = \left[\theta_{ub}, \theta_{ur} \right]$ are the reference gimbal angles given to the feedback linearization. 

\begin{figure*}
	\centering
		\includegraphics[width=0.7\textwidth]{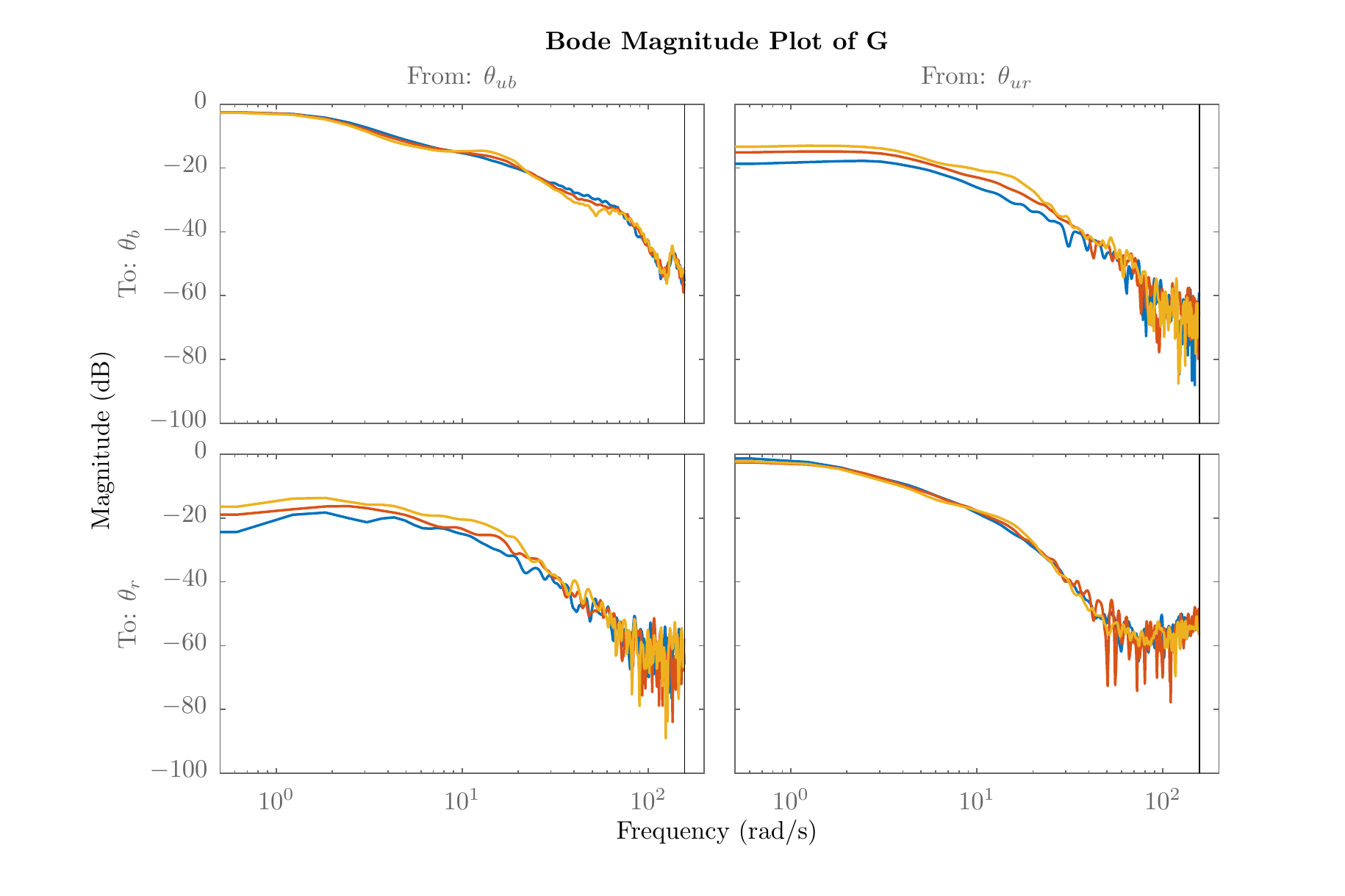}
	\caption{The identified frequency response of the blackbox model $G$ at different disk speeds. The blue line is the response at a disk speed of 300~rpm, red at 400~rpm and yellow at 500~rpm.}
	\label{fig:freqresp_gyro}
\end{figure*}

The inner loop is then taken as a black box model $G$ with 2 inputs and 2 outputs, and a single-channel identification is performed to calculate the frequency response of the new plant. A PRBS signal with an amplitude of $\pm10^{\circ}$, a length of 511 samples and a sampling time of 20~ms was applied for 4 periods to $\theta_{ub}$ and $\theta_{ur}$ respectively. The frequency response was calculated in Matlab using the $spa$ command with a Hann window length of 150. The identification was performed for the three different disk velocities $V= \left[300, 400, 500 \right] \text{rpm}$, resulting in three models $G = \left[G_1, G_2, G_3\right]$. The frequency responses are shown in Fig.~\ref{fig:freqresp_gyro}. It can be seen that the coupling and resonance modes become stronger at higher disk speeds.

\subsubsection{Control design formulation}

Based on the three frequency responses, a multivariable controller is designed. The goal is to decouple the system while also achieving good tracking performance of the reference angles $\theta^*$. Therefore, as objective function we choose to minimize the 2-norm $\| L-L_d  \|_2^2$ between the actual open-loop transfer function $L$ and desired open-loop transfer function $L_d=\frac{4}{s}I$,  
where a bandwidth of 4 rad/s is desired for the decoupled system.  The effect of the high frequency resonance mode is reduced by choosing a high pass filter for the complementary sensitivity function. To avoid input saturation a constant weighting is considered for the input sensitivity function $U=KS$.  The $H_\infty$ constraints are:
\begin{equation}
	 \|W_1T\|_{\infty} < 1 \quad ; \quad 
	\|W_2U\|_{\infty} < 1 	
\end{equation}
where $W_1(j\omega)=(0.2 j\omega+1)I$ and $W_2=0.05 I$. A 4th-order discrete-time controller with a sampling time of 0.04~s is chosen for this example. The controller includes an integrator, i.e. $F_y= (z - 1)I$. The matrix $Y$ is chosen to be diagonal. This choice of $Y_i$ greatly simplifies the calculation of the inverse and leads to the input channels having the same dynamics to every output. Note that the desired $L_d$ and the weighting filters can be in continuous-time, while the designed controller is in discrete-time. The fact that $W_1$ is not proper does not create any problem in practice because the constraints are evaluated for finite values of $\omega$.

The optimization problem is sampled using $N=500$ frequency points in the interval $\Omega_N = \left[ 10^{-1}, 25\pi \right]$ (the upper limit being the Nyquist frequency of the controller). The lower limit is chosen greater than zero in order to guarantee the boundedness of $L-L_d$. In fact a weighted two-norm of $L-L_d$ which is bounded is minimized.

The constraint sets are formulated for each of the three identified models $\left[G_1, G_2, G_3 \right]$, resulting in the following optimization problem : 
\begin{align*}
&\hspace{2cm}\min_{X,Y} \sum_{i = 1}^{3} \sum_{k = 1}^{N} \mbox{trace}[\Gamma_{k_i}] \\
&\text{subject to:}\\
&\left[\begin{array}{cc}Y^*Y_c+Y_c^*Y-Y_c^*Y_c & \quad(G_iX-L_dY)^* \\G_i X-L_dY & \Gamma_{k_i} \end{array}\right] (j \omega_k) > 0 \\
&\left[\begin{array}{cc} P_i^*P_{c_i} + P_{c_i}^*P_i - P_{c_i}^*P_{c_i} & (W_1G_iX)^* \\ W_1G_iX & I \end{array}\right] (j \omega_k) > 0 \\
&\left[\begin{array}{cc} P_i^*P_{c_i} + P_{c_i}^*P_i - P_{c_i}^*P_{c_i} & (W_2X)^* \\ W_2X & I \end{array}\right] (j \omega_k) > 0 \\
&k=1,\ldots,N \quad ; \quad i=1,2,3
\end{align*}

\begin{figure*}
	\centering
		\includegraphics[width=0.7\textwidth]{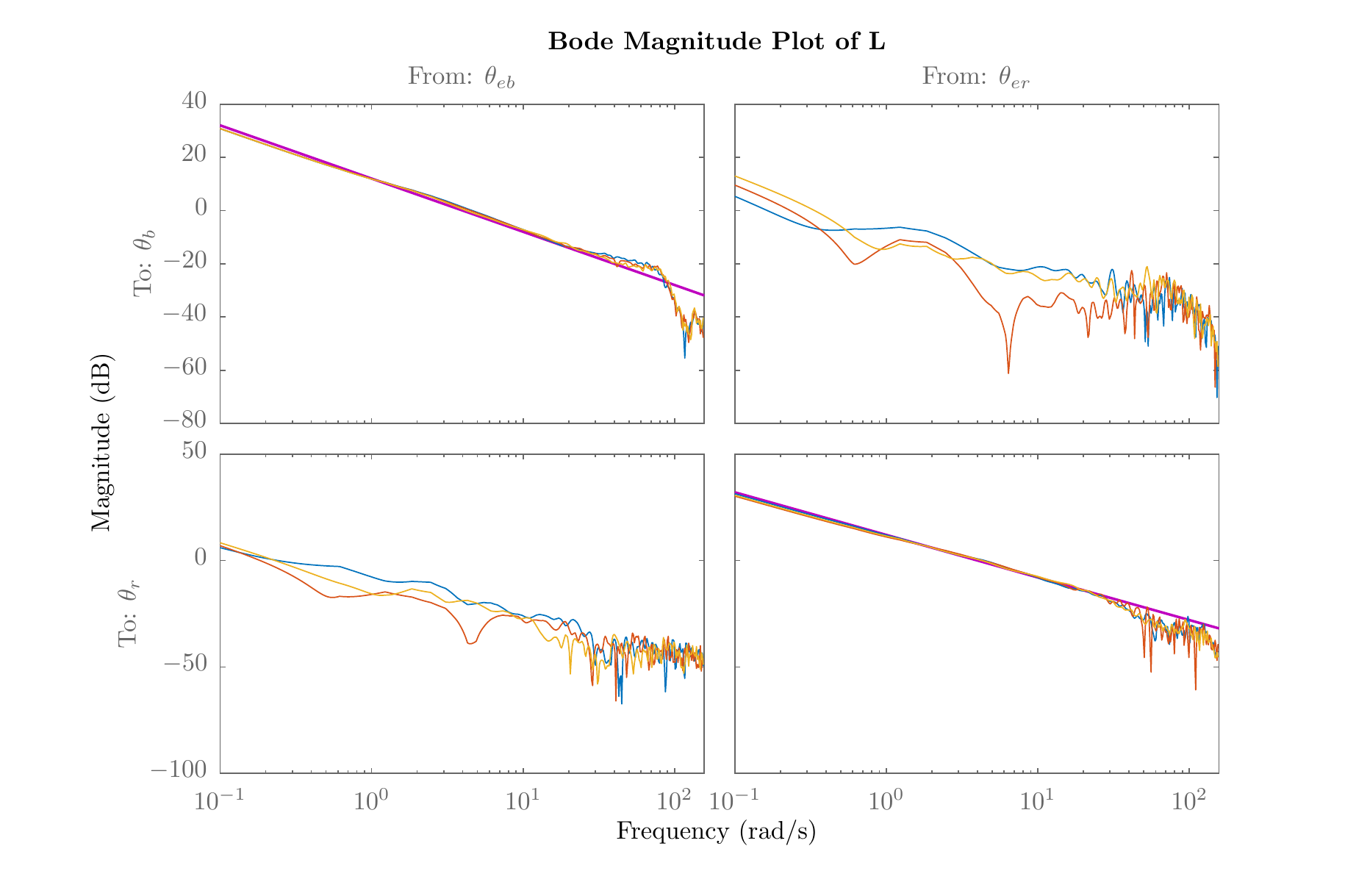}
	\caption{Bode magnitude plots of the open-loop transfer functions $L_d$ and $L_{1,2,3}$ for the three different plant models. The blue line is the actual response at a disk speed of 300~rpm, red at 400~rpm and yellow at 500~rpm. The desired $L_d$ is shown in purple.}
	\label{fig:gyro_bodemag_L_Ld}
\end{figure*}

As the gyroscope is a stable system, the initial controller was chosen by setting the poles of the controller to 0 and choosing a small enough gain: 
\begin{equation}
	X_c= 0.01 I \quad ; \quad Y_c=z^4(z-1)I 
\end{equation}
The iteration converges to a final controller in 10 steps. The bode magnitude plots of $L_d$ and $L_{1,2,3}$ for the three different plant models are shown in Fig.~\ref{fig:gyro_bodemag_L_Ld}. It can be seen that the designed controller approximates the desired loop shape well in low frequencies, and that the coupling has been reduced.

\subsubsection{Experimental results}

\begin{figure}
	\centering
		\includegraphics[width=0.9\columnwidth]{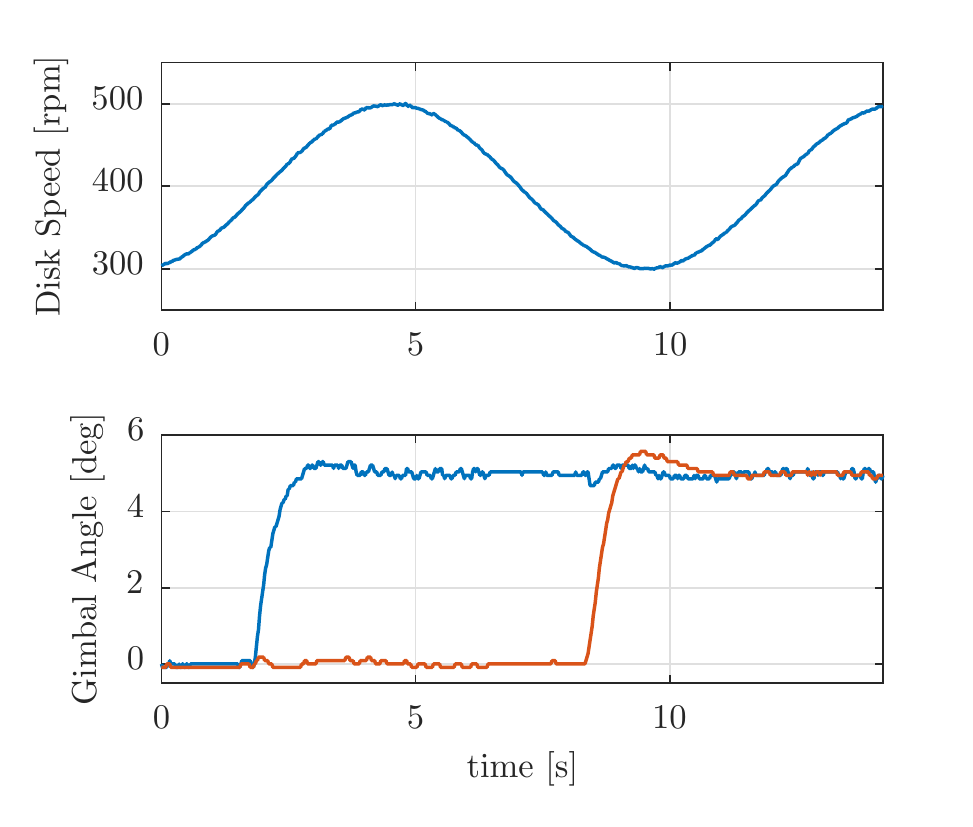}
	\caption{Step response of the blue and red gimbal angles during a varying disk velocity.}
	\label{fig:gyro_steps}
\end{figure}

To validate the results, the controller was implemented in Labview and applied to the experimental setup. The step responses of the blue and red gimbal angle were measured for varying disk speeds, and the results are shown in Fig.~\ref{fig:gyro_steps}. It can be seen that the decoupling is good, and that the multimodel uncertainty introduced by the varying disk speed is handled well. The rise time is 0.625~s for the blue and 0.486~s for the red gimbal angle, which is close to the desired rise time of 0.55~s. A slight overshoot can be observed especially for the red gimbal angle, which is likely due to the nonlinearities present in the system.

\section{Conclusions}\label{sec7}
The frequency response of a multivariable system can be obtained through several experiments. This data can be used directly to compute a high performance controller without a parametric identification step. The main advantage is that there will be no unmodeled dynamics and that the uncertainty originating from measurement noise can be straightforwardly modeled through the weighting frequency functions. A unified convex approximation is used to convexify the $H_\infty$, $H_2$ and loop shaping control problems. Similar to the model-based approaches, this convex approximation relies on an initial stabilizing controller. Several initialization techniques are discussed and an iterative algorithm is proposed that converges to a local optimum of the original non-convex problem. Compared to the other frequency-domain data-driven approaches, the proposed method has a full controller parametrization and also covers $H_2$ and loop shaping control design with a new closed-loop stability proof. 

%
\bibliography{linear}
\bibliographystyle{automatica}

\noindent {\bf Alireza Karimi} received his PhD in 1997 from Institut National Polytechnique de Grenoble (INPG) in France. He was Assistant Professor at Electrical Engineering Department of Sharif University of Technology in Teheran from 1998 to 2000. He is currently Senior Scientist at the Automatic Control Laboratory of Ecole Polytechnique F\'ed\'erale de Lausanne (EPFL), Switzerland. He was an Associate Editor of European Journal of Control from 2004 to 2013. His research interests include closed-loop identification, data-driven controller tuning approaches and robust control. 

\noindent{\bf Christoph Kammer} is a doctoral student in electrical engineering at the Automatic Control Laboratory at EPFL. He received his master's degree in mechanical engineering from ETHZ in 2013. His main interests lie in multivariable robust control and its applications in power grids.

\end{document}